\documentclass[12pt
]{article}
\usepackage{amsmath,amsfonts,latexsym,mathrsfs, amsthm,amssymb}
\makeatletter
  
  \@addtoreset{equation}{section}
\makeatother
\newtheorem{proposition}{Proposition}[section]
\newtheorem{theorem}{Theorem}[section]

\newtheorem{remark}{Remark}[section]
\newtheorem{corollary}{Corollary}[section]

\date{ }

\begin{document}

\title{A remark on the Lagrangian formulation of 
optimal transport with a non--convex cost
}

\author{
Toshio Mikami\thanks{Partially supported by JSPS KAKENHI Grant Number 19K03548.
}\qquad Haruka Yamamoto
\\
Department of Mathematics, Tsuda University
}

\maketitle

\begin{abstract}
We study the Lagrangian formulation of  a class of the 
Monge--Kantorovich optimal transportation problem.
It can be considered a stochastic optimal transportation problem for 
absolutely continuous stochastic processes.
A cost function and stochastic processes under consideration 
is not convex and  have essentially bounded time derivatives almost surely, respectively.
This paper is a continuation of the second author's master thesis.
\end{abstract} 

Keywords:  Lagrangian formulation, non--convex cost, Monge--Kantorovich problem, stochastic optimal transport

AMS subject classifications: 49Q22, 93E20, 49K45

\section{Introduction}
\label{intro}
For $d\ge 1$, let
$\mathcal{P}(\mathbb{R}^d)$ denote the space of all Borel probability measures on $\mathbb{R}^d$ endowed with weak topology.
For $P_0, P_1\in \mathcal{P}(\mathbb{R}^d)$,
let $\Pi(P_0, P_1)$ denote the set of $\mathbb{R}^d\times \mathbb{R}^d$--valued random variables
$(X_0,X_1)$ defined on a possibly different probability space such that
$P^{X_i}=P_i, i=0,1$.
Here $P^X$ denotes the probability distribution of a random variable $X$.
Throughout this paper, the probability space under consideration is not fixed.

Let $L:\mathbb{R}^d\rightarrow [0,\infty)$ be Borel measurable.
The following is a typical  Monge--Kantorovich optimal transportation problem:
for  $P_0, P_1\in \mathcal{P}(\mathbb{R}^d)$,
\begin{equation}\label{1.1}
T(P_0,P_1):=\inf\{E[L(X_1-X_0)]:(X_0,X_1)\in \Pi(P_0, P_1)\}
\end{equation}
(see, e.g. \cite{RR, V} and the references therein).

In the case where $L(u)=|u|^p$ for $p>0$, we denote (\ref{1.1}) by $T_{p}(P_0,P_1)$.
In the case where $d=1$, the minimizer of $T_{p}, p\ge 1$ was obtained in \cite{D}.
In the case where $d\ge 2$, the minimizer of $T_{2}$ was obtained in \cite{B1, B2} and the generalization to more general costs, including concave ones, was given in \cite{GM}.
The probabilistic proof of the existence and the uniqueness of the minimizer of $T_{2}$ via
the stochastic control approach was given in \cite{MMonge}
by the zero--noise limit of Schr\"odinger's problem
 (see \cite{S. B,S1, S2} for Schr\"odinger's problem and \cite{Jamison1975,L1,L2, M2021,Zambrini} for related topics).
Schr\"odinger's problem is also called the entropic regularized optimal transport these days and plays a crucial role in data science (see 
\cite{CGT,Cuturi,PC} and the reference therein).
Notice that $\mathbb{R}^d\ni u\mapsto |u|^p$ is convex for $p\in [1,\infty)$ and 
is not concave for $p\in (0,1)$.

Let $\mathcal{A}(P_0, P_1)$
 denote the set of stochastic processes $X(\cdot)$ such that 
 $$X(\cdot)\in AC([0,1]), {\rm a.s.},\quad (X(0),X(1))\in \Pi(P_0, P_1).$$
Here, for $T>0$, $AC([0,T])$ denotes the space of all absolutely continuous functions from $[0,T]$ to $\mathbb{R}^d$.
We consider the following stochastic optimal transport for absolutely continuous stochastic processes:
\begin{equation}\label{1.1.0}
V(P_0,P_1):=\inf\left\{E\left[\int_0^1 L(X'(t))dt\right]: X(\cdot) \in 
\mathbf{\mathcal{A}}(P_0, P_1)\right\},
\end{equation}
where $\displaystyle X'(t):=\frac{d}{dt}X(t)$.

We briefly describe the relation between $T$ and $V$.
The following holds without any assumption:
\begin{equation}\label{1.3}
V(P_0,P_1)\le T(P_0,P_1).
\end{equation}
If $L:\mathbb{R}^d\rightarrow [0,\infty)$ is convex,
then the equality holds in (\ref{1.3}) and 
$V(P_0,P_1)$ can be considered the Lagrangian formulation of $T(P_0,P_1)$
(see, e.g. \cite{M2021} for the proof of (\ref{1.3})
and also \cite{Mikami2002ECP} for related topics).\\
If (i) $L(ru)\ge rL(u), 0<r<1, u\in \mathbb{R}^d$; (ii) $L(u)/|u|\to 0, |u|\to\infty$, and (iii) $T(P_0,P_1)$ is finite,
then
\begin{equation}\label{1.6}
V(P_0,P_1)=0
\end{equation}
(see Appendix for the proof and also Theorem \ref{thm2.1} and Corollary \ref{co2.2.0617} in section 2).
A typical example of $L$ is $|u|^p, p\in (0,1)$
(see Remark \ref{rk2.1.0} in section \ref{section3} for more examples).

(\ref{1.3})--(\ref{1.6}) imply that, 
to study the Lagrangian formulation of $T(P_0,P_1)$ when  $L$ is not convex,
we have to modify a cost function or restrict a class of stochastic processes in (\ref{1.1.0}).

We first modify a cost function and give two Lagrangian formulations for $T(P_0,P_1)$.
For $t>0, \varphi\in L^\infty([0,t])$,
\begin{equation}\label{1.8.0416}
1\le N_1(\varphi)_t:=
\begin{cases}
\displaystyle \frac{t||\varphi||_{\infty,t}}{|\int_0^t\varphi(s)ds|},
&\hbox{if }\int_0^t\varphi(s)ds\ne 0,\\
1, &\hbox{otherwise},
\end{cases}
\end{equation}
\begin{equation}\label{2.4.0427}
1\le N_2(\varphi)_t:=
\begin{cases}
\displaystyle \frac{t||\varphi||_{\infty,t}}{||\varphi||_{1,t}},
&\hbox{if }||\varphi||_{1,t}>0,\\
1, &\hbox{otherwise},
\end{cases}
\end{equation}
where 
$$||\varphi||_{\infty,t}:=ess.sup \left\{|\varphi (s)|:0\le s\le t\right\},
\quad ||\varphi||_{1,t}:=\int_0^t |\varphi(s)|ds.$$
For simplicity, $||\varphi||_{p}:=||\varphi||_{p,1}$ for $p=1,\infty$ and 
$N_i(\varphi):=N_i(\varphi)_1$, $i=1,2$.

Let $\ell:[0,\infty)\rightarrow [0,\infty)$ and $\ell (0)=0$ (see (A1) in section \ref{section3}).
For $i=1,2$,
\begin{equation}
L_i(t,\varphi):=
N_i(\varphi)\ell \left( \frac{|\varphi(t)|}{N_i(\varphi)}\right),\quad (t,\varphi)\in
[0,1]\times L^\infty([0,1]),
\end{equation}
\begin{equation}\label{1.10.0416}
\tilde V_i(P_0,P_1) :=\inf\left\{E\left[\int_0^1 L_i(t,X')dt\right]:
  X(\cdot) \in \mathcal{A}_\infty(P_0,P_1)\right\},
\end{equation}
where 
$$\mathcal{A}_\infty(P_0,P_1):=\{X(\cdot) \in \mathcal{A}(P_0,P_1):
||X'||_\infty<\infty  \hbox{  a.s.}\}.$$
We show that 
$$T(P_0,P_1)=\tilde V_i(P_0,P_1),\quad i=1,2,$$
under different assumptions  (see Theorems \ref{thm2.1.0416}--\ref{co2.1.0427} in section \ref{section3}).

\begin{remark}
If $L:\mathbb{R}^d\rightarrow [0,\infty)$ is convex, $L(0)=0$ and $V(P_0,P_1)$ is finite, then the following holds (see Appendix for the proof): for $i=1,2$,
\begin{eqnarray}\label{1.9.0416}
&&V(P_0,P_1)=T(P_0,P_1)\\
&=&\inf\left\{E\left[\int_0^1 \frac{1}{N_i(X')}L\left(N_i(X')X'(t)\right)dt\right]: X(\cdot) \in \mathcal{A}_\infty(P_0,P_1)\right\}.\notag
\end{eqnarray}
\end{remark}

For  $ X(\cdot) \in \mathcal{A}_\infty(P_0,P_1)$ and $t\ge 0$,
$L_i(t,X')\ge 0$ and $=0$ if $X'(t)=0$.
When we consider minimizers of $\tilde V_i$,
we assume that $\ell(u)>0, u>0$ so 
that we only have to consider $X(\cdot) $ such that $X(t)=X(0)$ if and only if $||X'||_{1,t}=0$.
In particular, we can assume that the following holds:
\begin{equation}\label{2.5.0610}
N_1(X')\ge N_2(X')\ge 1,\quad {\rm a.s.},
\end{equation}
which implies the following: 
\begin{equation}\label{1.12.0608}
|X'(t)|\ge \frac{|X'(t)|}{N_2(X')}\ge \frac{|X'(t)|}{N_1(X')},\quad L_1(t,X')\ge L_2(t,X'),
\end{equation}
provided 
$\ell (ru)\ge r\ell (u)$ for $(r,u)\in (0,1)\times (0,\infty)$ (see (\ref{3.13.0428})
and also Theorems \ref{thm2.1.0416}--\ref{co2.1.0427} and Proposition \ref{pp2.1} in section \ref{section3}).
The following also holds (see Appendix for the proof):
\begin{eqnarray}\label{1.11.0508}
&&\tilde V_i(P_0,P_1)\\
&=&\inf\left\{E\left[\int_0^\tau\ell \left(|X'(t)|\right)dt\right]\right.: \tau=\tau(\omega)\ge 1,X(\cdot) \in AC([0,\tau]),\notag\\
&&\qquad \left. 
N_i(X')_\tau=\tau, {\rm a.s.}, (X(0),X(\tau) )\in \Pi(P_0, P_1)\right\}.\notag
\end{eqnarray}

Next, we consider a restricted class of absolutely continuous stochastic processes with almost surely essentially bounded time derivatives.
For $ P_0, P_1\in \mathcal{P}(\mathbb{R}^d)$, and $B\subset \mathcal{P}([0,\infty))$,
\begin{eqnarray*}
&&\mathcal{A}_\infty(P_0, P_1;B)\\
&:=&
 \{(X(\cdot),M): X(\cdot)\in \mathcal{A}_\infty(P_0, P_1) ,P^{M}\in B,  ||X'||_\infty\le M, {\rm a.s.}\},\\
 &&\Pi_\infty (P_0, P_1;B)\\
&:=&
 \{(X_0,X_1,M): (X_0,X_1)\in \Pi(P_0, P_1), P^{M}\in B,  |X_1-X_0|\le M, {\rm a.s.}\}.
 \end{eqnarray*}
\begin{eqnarray}
&&V(P_0,P_1;B)\label{1.9}\\
&:=&\inf\left\{E\left[\int_0^1\ell(|X'(t)|)dt\right]: (X(\cdot),M)\in 
\mathcal{A}_\infty(P_0, P_1;B)\right\},\notag\\
&&T^V(P_0,P_1;B)\\
&:=&\inf\{E[\ell (M)M^{-1}|X_1-X_0|;M>0]:
(X_0,X_1,M)\in \Pi_\infty (P_0, P_1;B)\}.\notag
\end{eqnarray}
We show that the following holds (see Theorem \ref{thm2.1} in section \ref{section3}):
$$V(P_0,P_1;B)=T^V(P_0,P_1;B).$$ 
It is a continuation of the second author's master thesis \cite{Yamamoto}
in which she only considered the case where $L(u)=|u|^p, p\in (0,1)$ and  $B$ is a set of a delta measure.

A generalization of our result to the case where stochastic processes under consideration 
are semimartingales is the first step to the theory of  stochastic optimal transport with a non--convex cost and is our future project.

We state our results in section \ref{section3} and prove them in section \ref{section4}. 
In Appendix, we give  the proofs for (\ref{1.6}),  (\ref{1.9.0416}), and (\ref{1.11.0508}) for the sake of completeness.

\section{Main result}\label{section3}

In this section, we state our results.
We first state the assumptions.

\noindent
(A1). (i) 
$\ell:[0,\infty)\rightarrow [0,\infty)$, $\ell (0)=0$,
\begin{equation}\label{a2.1}
\ell (ru)\ge r\ell(u),\quad (r,u)\in (0,1)\times (0,\infty).
\end{equation}
(ii)  In (\ref{a2.1}), the equality does not hold for any $(r,u)\in (0,1)\times (0,\infty)$.
(iii) $\ell (u)>0, u>0$.

\noindent
(A2). (i) $\ell:[0,\infty)\rightarrow [0,\infty)$ is non--decreasing.
(ii) $\ell:[0,\infty)\rightarrow [0,\infty)$ is strictly increasing.
(iii) $\ell\in C([0,\infty))$ and $\ell(u)\to\infty$, as $|u|\to\infty$.

\noindent

We state remarks on (A1)--(A2).

\begin{remark}\label{rk2.1.0}
(i) (\ref{rk2.1.0}) and (A1,ii) mean that $(0,\infty)\ni u\mapsto \ell (u)/u$ is non--increasing
and is strictly decreasing, respectively.
In particular, (A1,ii) implies (A1,iii), provided $\ell (u)\ge 0$.\\
(ii) If $\ell$ is concave and $\ell (0)=0$, then (\ref{a2.1}) holds.
If $\ell$ is strictly convex and $\ell (0)=0$, then (\ref{a2.1}) does not hold.
$\ell (u)=u$ satisfies (A1,i), but not  (A1,ii).\\
(iii) 
$$\ell (u)=
\begin{cases}
2u\exp (-u), &0\le u< 1,\\
u\exp (-u), &u\ge 1
\end{cases}
$$
is concave on $[0,1)$ and $[1,2]$ and is convex on $[2,\infty)$.
It is strictly increasing and strictly decreasing on $[0,1)$ and $[1,\infty)$, respectively.
It is not continuous at $u=1$ and satisfies (A1).\\
(iv) (\ref{a2.1}) and (A2,i) imply that $\ell\in C((0,\infty))$ since, if $0<h<u$, 
$$\frac{\ell(u)}{u+h}\le \frac{\ell(u+h)}{u+h}\le \frac{\ell(u)}{u}\le \frac{\ell(u-h)}{u-h}\le \frac{\ell(u)}{u-h}.$$

\end{remark}

We describe a list of notations of the sets of minimizers.

$\Pi_{T,opt}(P_0,P_1):=$the set of  minimizers of $T(P_0,P_1)$.

$\Pi_{T^V,opt}(P_0,P_1;m):=$the set of minimizers of $T^V(P_0,P_1;m)$.

$\mathcal{A}_{i,opt}(P_0,P_1):=$the set of minimizers of $\tilde V_i(P_0,P_1), i=1,2$.

$\mathcal{A}_{opt}(P_0,P_1;m):=$the set of minimizers of $V(P_0,P_1;m)$.


\bigskip

We say that  $A\subset [0,1]$ is a random measurable set 
if and only if there exists a $\{0,1\}$--valued stochastic process $\{\eta(t,\omega)\}_{0\le t\le 1}$ defined on a probability space $(\Omega ,\mathcal{F}, P)$ such that 
$$[0,1]\times \Omega\ni (t,\omega)\mapsto \eta(t,\omega)\in \{0,1\}$$
is jointly measurable and $A=A(\omega) =\eta(\cdot,\omega)^{-1}(1)$, i.e.
$\eta(t,\omega)=I_{A(\omega)}(t)$, where $I_B(x)=1, x\in B; =0,x\not\in B$.
It is easy to see that the Lebesgue measure $|A(\omega) |=\int_0^1 I_{\{1\}}(\eta(t,\omega))dt$ is a random variable.

For $x,y\in \mathbb{R}^d$, a Lebesgue measurable set $A\subset [0,1]$,
and $t\in [0,1]$,
\begin{equation}\label{2.3.0525}
X(t;x,y,A):=
\begin{cases}
\displaystyle x+\frac{|A\cap [0,t]|}{|A|}(y-x),& \hbox{ if } x\ne y,|A|>0,\\
x,&\hbox{otherwise}.
\end{cases}
\end{equation}

The following gives the relation between $T(P_0,P_1)$ and $\tilde V_1(P_0,P_1)$ 
(see (\ref{1.1}) and (\ref{1.10.0416}) for notation).

\begin{theorem}\label{thm2.1.0416}
Suppose that  (A1,i) holds.
Then for any $P_0, P_1\in \mathcal{P}(\mathbb{R}^d)$, the following holds.\\
(i)
\begin{equation}\label{2.1.0415}
T(P_0,P_1)=\tilde V_1(P_0,P_1).
\end{equation}
(ii) If $(X_0,X_1)\in \Pi_{T,opt}(P_0,P_1)$ and a random measurable set $A\subset [0,1]$
are defined on the same probability space and if 
$$P(|A|>0|X_0\ne X_1)=1,$$
then
$X(\cdot;X_0,X_1,A)\in \mathcal{A}_{1,opt}(P_0,P_1)$.


\noindent
(iii) If $X(\cdot)\in \mathcal{A}_{1,opt}(P_0,P_1)$,  then $(X(0),X(1))\in \Pi_{T,opt}(P_0,P_1)$.
Suppose, in addition, that (A1,ii) holds.  Then  
$X(\cdot)=X(\cdot;X(0),X(1),(X')^{-1}(\mathbb{R}^d\backslash\{0\}))$,
where 
$$(X')^{-1}(\mathbb{R}^d\backslash\{0\}):=\{t\in [0,1]:X'(t)\ne 0\}.$$
\end{theorem}

\begin{remark}
In the case where  $L$ is strictly convex, for an optimal path $X(\cdot)$ of $V(P_0,P_1)$,
$X(\cdot)=X(\cdot;X(0),X(1),[0,1])$ by Jensen's inequality (see, e.g. \cite{M2021}).
In particular, it moves at constant velocity.
Theorem \ref{thm2.1.0416} implies that under  (A1), 
an optimal path $X(\cdot)\in \mathcal{A}_{1,opt}(P_0,P_1)$ can stop 
even randomly.
But when it moves, the velocity is constant in $t$ and  can be random.
\end{remark}

Under (A1,i,iii), $\tilde V_1\ge\tilde V_2$  (see (\ref{3.13.0428})).
The following implies that equality holds under an additional assumption (A2,i).

\begin{theorem}\label{co2.1.0427}
Suppose that   (A1,i,iii) and (A2,i) hold.
Then for any $P_0, P_1\in \mathcal{P}(\mathbb{R}^d)$, the following holds.\\
(i)
\begin{equation}\label{2.6.0427}
\tilde V_1(P_0,P_1)=\tilde V_2(P_0,P_1).
\end{equation}
(ii) 
Suppose, in addition, that 
(A2,ii) holds. Then  
$$\mathcal{A}_{1,opt}(P_0,P_1)=\mathcal{A}_{2,opt}(P_0,P_1).$$
In particular, for any $X(\cdot)\in \mathcal{A}_{1,opt}(P_0,P_1)$,
\begin{equation}\label{2.7.0427}
N_1(X')=N_2(X').
\end{equation}
\end{theorem}

The following implies that  Theorem \ref{co2.1.0427} does not necessarily hold without (A2,i)
(see Remark \ref{rk2.1.0}, (iii) for an example and also Theorem \ref{thm2.1.0416}). 

\begin{proposition}\label{pp2.1}
Suppose that there exists $r_0>0$ such that $\ell$ is strictly decreasing on $[r_0,\infty)$.
Then for any $P_0, P_1\in \mathcal{P}(\mathbb{R}^d)$
for which $T(P_0,P_1)$ has a minimizer $(X_0,X_1)$ such that 
$P(|X_1-X_0|\ge r_0)>0$, the following holds:
\begin{equation}
T(P_0,P_1)>\tilde V_2(P_0,P_1).
\end{equation}
\end{proposition}

For $f\in C_b (\mathbb{R}^d)$,
\begin{equation}
f^\ell(x):=\inf\{\ell (|y-x|)+f(y)|y\in \mathbb{R}^d\},\quad x\in \mathbb{R}^d.
\end{equation}
From (i) in Theorems \ref{thm2.1.0416}--\ref{co2.1.0427}, we easily obtain the following  and omit the proof (see the proof of Theorem 2.1 in \cite{2-MikamiSImple}).

\begin{corollary}\label{co2.1.0611}
Suppose that  (A2,iii) holds.
Suppose also that
``(A1,i)'' or  ``(A1,i,iii) and (A2,i)'' 
hold.
Then for  $i=1$ or  $2$, the following holds, respectively:
for any $P_0\in \mathcal{P}(\mathbb{R}^d)$ such that $P_0(dx)\ll dx$ and any $f\in C_b (\mathbb{R}^d)$,
\begin{eqnarray}\label{2.6.0609} 
&&\inf\left\{E\left[\int_0^1L_i(t,X)dt+f(X(1))\right]:X\in \mathcal{A}_\infty(P_0,P^{X(1)})\right\}\qquad\\
&=&\int_{\mathbb{R}^d}f^\ell(x)P_0(dx).\notag
\end{eqnarray}
\end{corollary}

\begin{remark}
(\ref{2.6.0609}) is a finite--time horizon optimal control problem for absolutely continuous stochastic processes (see \cite{2-FS} for stochastic control theory) and
the l. h. s. can be also written as follows:
$$\inf\left\{\tilde V_i(P_0,P)+\int_{\mathbb{R}^d}f(x)P(dx):P\in\mathcal {P}(\mathbb{R}^d)\right\}.$$
\end{remark}

The following gives the relation between $V(P_0,P_1;B)$ and $T^V(P_0,P_1;B)$.

\begin{theorem}\label{thm2.1}
Suppose that (A1,i) holds.
Then for any $P_0, P_1\in \mathcal{P}(\mathbb{R}^d)$ and $B\subset\mathcal{P}([0,\infty))$, the following holds.\\
(i)
\begin{equation}\label{2.1}
V(P_0,P_1;B)=T^V(P_0,P_1;B).
\end{equation}
(ii) If $(X_0,X_1,M)\in \Pi_{T^V,opt}(P_0,P_1;B)$
and a random measurable set $A\subset [0,1]$ are defined on the same probability space
and if 
$$P\left(|A|=\frac{|X_1-X_0|}{M}\biggl|M>0\right)=1,$$
then
$(X(\cdot;X_0,X_1,A),M)\in \mathcal{A}_{opt}(P_0,P_1;B)$.

\noindent
(iii) If $(X(\cdot),M)\in \mathcal{A}_{opt}(P_0,P_1;B)$, then
$(X(0),X(1),M)\in \Pi_{T^V,opt}(P_0,P_1;B)$.
Suppose, in addition, that (A1,ii) holds.  
Then
\begin{equation}\label{2.2}
P\left(|(X')^{-1}(\mathbb{R}^d\backslash\{0\})|=\frac{|X(1)-X(0)|}{M}\biggl|M>0\right)=1,
\end{equation} 
and  $X(\cdot)=X(\cdot;X(0),X(1),(X')^{-1}(\mathbb{R}^d\backslash\{0\}))$.
\end{theorem}

\begin{remark}\label{rk1.1}
Even if $\ell$ is lower semicontinuous,
$$AC([0,1])\ni \varphi\mapsto \int_0^1 \ell (|\varphi'(t)|)dt$$
is not necessarily lower semicontinuous 
in the  supnorm.
In particular, it is not trivial if $\mathcal{A}_{opt}(P_0,P_1;B)$ is not empty.
\end{remark}

(A1,i) implies that $\ell(u)/u$ is convergent 
 as $u\to\infty$ (see Remark \ref{rk2.1.0}, (i)):
$$C_\ell:=\lim_{u\to\infty}\frac{\ell (u)}{u}.$$
In particular, the following holds from  Theorem \ref{thm2.1}, (i).

\begin{corollary}\label{thm2.5.0628}
Suppose that (A1,i) holds.
Then for any $P_0, P_1\in \mathcal{P}(\mathbb{R}^d)$, the following holds:
\begin{equation}
\inf\left\{E\left[\int_0^1 \ell (|X'(t)|)dt\right]:
X(\cdot) \in  \mathcal{A}_\infty(P_0, P_1)\right\}=C_\ell \cdot T_1(P_0,P_1).\label{2.17.0628}
\end{equation}
In particular, if $P_0\ne P_1$ and the l. h. s. of (\ref{2.17.0628}) has a minimizer, then  
\begin{equation}\label{2.13.0718}
\inf\left\{u>0:\frac{\ell (u)}{u}=C_\ell\right\}<\infty.
\end{equation}
\end{corollary}

$r\mapsto V(P_0,P_1;\{\delta_r\})$ is non--increasing and converges to $V(P_0,P_1;\{\delta_M\}_{M>0})$,
as $r\to\infty$, where $\delta_r$ denotes the delta measure on $\{r\}$.
In particular, we easily obtain the following from  Theorem \ref{thm2.1}, (i) and we omit the proof.

\begin{corollary}\label{co2.2.0617}
Suppose that (A1,i) holds.
Then for any $P_0, P_1\in \mathcal{P}(\mathbb{R}^d)$ with bounded supports, the following holds:
for any $r\ge \sup\{|x_0-x_1|;x_i\in supp(P_i),i=0,1\}$, 
\begin{equation}
V(P_0,P_1;\{\delta_r\})=\frac{\ell (r)}{r}T_1(P_0,P_1).
\end{equation}
In particular, 
\begin{equation}\label{2.16.0617}
V(P_0,P_1;\{\delta_M\}_{M>0})
=C_\ell \cdot T_1(P_0,P_1)
\end{equation}
and the left--hand sides of  (\ref{2.17.0628})   and  (\ref{2.16.0617})  coincide.
\end{corollary}

\begin{remark}
For $a\ge 0$, $\ell(u)=au+1-\exp (-u)$ is concave, satisfies (A1,i), and  $C_\ell =a$.
\end{remark}


\section{Proofs of results in section 2}\label{section4}

In this section, we prove our results.
When it is not confusing, we omit ``${\rm a.s.}$'' for the sake of simplicity.

We first prove Theorem \ref{thm2.1.0416}.
\begin{proof}
(Theorem \ref{thm2.1.0416})
We first prove (i).
We prove
\begin{equation}\label{3.1.0416}
T(P_0,P_1)\le \tilde V_1(P_0,P_1).
\end{equation}
Suppose that $X(\cdot)\in \mathcal{A}_\infty(P_0,P_1)$.
If $X(1)\ne X(0)$, then 
\begin{equation}\label{3.3.0416}
N_1\int_0^1\ell \left(\frac{|X' (t)|}{N_1}\right)dt \ge \ell\left(|X(1)-X(0)|\right),
\end{equation}
where $N_1=N_1(X')$ (see (\ref{1.8.0416}) for notation).
Indeed, from (A1,i),
\begin{eqnarray}\label{3.2.0416}
\ell \left(\frac{|X' (t)|}{N_1}\right)
&=&\ell\left(\frac{|X' (t)|}{N_1|X(1)-X(0)|}|X(1)-X(0)|\right)\\
&\ge &\frac{|X' (t)|}{N_1|X(1)-X(0)|}\ell\left(|X(1)-X(0)|\right),\quad dt-{\rm a.e.}\notag
\end{eqnarray}
since
$$\frac{|X'(t)|}{N_1|X(1)-X(0)|}=\frac{|X'(t)|}{||X'||_\infty}\le 1,\quad dt-{\rm a.e.}.$$
Besides,
\begin{equation}\label{3.4.0416}
||X'||_1 \ge |X(1)-X(0)|.
\end{equation}
If $X(1)=X(0)$, then (\ref{3.3.0416}) holds trivially.
(\ref{3.3.0416}) implies (\ref{3.1.0416}) immediately.

We prove  
\begin{equation}\label{3.5.0416}
\tilde V_1(P_0,P_1)\le T(P_0,P_1).
\end{equation}
Suppose that $(X_0,X_1)\in \Pi(P_0, P_1)$.
\begin{equation}\label{3.6.0416}
X(t):=X_0+t(X_1-X_0),\quad 0\le t\le 1.
\end{equation}
Then $X(\cdot)\in \mathcal{A}_\infty (P_0, P_1), N_1=N_1(X')=1$, and
\begin{equation}\label{3.7.0416}
\ell( |X_1-X_0|)=N_1\int_0^1  \ell \left(\frac{|X' (t)|}{N_1}\right)dt,
\end{equation}
which implies (\ref{3.5.0416}).

We prove (ii).
We write $X(\cdot;X_0,X_1,A)=X(\cdot;A)$ for simplicity.
$X(\cdot;A)\in \mathcal{A}_\infty (P_0, P_1)$ and (\ref{3.7.0416}) with $X(\cdot)=X(\cdot;A)$ holds since
$$|X'(t;A)|=I_A(t)\frac{|X_1-X_0|}{|A|}, \quad dt-{\rm a.e.,}\quad
N_1(X'(\cdot ;A))=\frac{1}{|A|},
$$
provided $X_1\ne X_0, |A|>0$.
 (\ref{2.1.0415}) and (\ref{3.7.0416})  with $X(\cdot)=X(\cdot;A)$ imply (ii). 
 
 (\ref{2.1.0415}) and (\ref{3.3.0416}) imply  the first part of (iii).
 We prove that $X(\cdot)=X(\cdot;X(0),X(1),(X')^{-1}(\mathbb{R}^d\backslash\{0\}))$.
For $X(\cdot)\in \mathcal{A}_{1,opt}(P_0,P_1)$, 
if $X (1)\ne X(0)$, then the equality holds in (\ref{3.3.0416})--(\ref{3.4.0416}).
In particular, the following holds:
\begin{eqnarray}\label{3.8.0416}
||X'||_1 &=&|X (1)-X(0)|,\quad {\rm a.s.,}\\
|X' (t)|&=&0, N_1|X (1)-X(0)|,\quad dtdP{\rm -a.e.}
\notag
\end{eqnarray}
from (A1,ii), where $N_1:=N_1(X')$.
Notice that the equalities in (\ref{3.8.0416}) hold if $X (1)=X(0)$ (see an explanation above (\ref{2.5.0610})
and Remark \ref{rk2.1.0}, (i)).
The following completes the proof:
\begin{equation}\label{3.9.0416}
X' (t)=\frac{X (1)-X(0)}{ |(X')^{-1}(\mathbb{R}^d\backslash\{0\})|}, \quad \hbox{ on }(X')^{-1}(\mathbb{R}^d\backslash\{0\}),  \quad dtdP{\rm - a.e.}.
\end{equation}
We prove (\ref{3.9.0416}).
(\ref{3.8.0416}) implies that for $P$-- almost all $\omega$, 
there exists $Z=Z(\omega) $ such that $|Z|=1$ and 
\begin{equation}\label{3.10.0416}
X' (t)=N_1|X (1)-X(0)|Z, \quad \hbox{ on }(X')^{-1}(\mathbb{R}^d\backslash\{0\}),  \quad dt{\rm -a.e.}.
\end{equation}
In particular,
\begin{equation}\label{3.11.0416}
X (1)-X(0)=N_1|X (1)-X(0)|Z\times |(X')^{-1}(\mathbb{R}^d\backslash\{0\})|.
\end{equation}
(\ref{3.10.0416})--(\ref{3.11.0416}) imply  (\ref{3.9.0416}).
\end{proof}

We prove Theorem \ref{co2.1.0427}.
\begin{proof}
(Theorem \ref{co2.1.0427})
We first prove (i).
We prove
\begin{equation}\label{3.12.0427}
\tilde V_2(P_0,P_1)\le \tilde V_1(P_0,P_1).
\end{equation}
For $X(\cdot)\in \mathcal{A}_\infty(P_0,P_1)$
such that $||X'||_\infty =0$ if $X(1)=X(0)$ and hence $N_2\le  N_1$,
from (A1,i),
\begin{equation}\label{3.13.0428}
N_1\ell \left(\frac{|X' (t)|}{N_1}\right)
=N_1\ell \left(\frac{|X' (t)|}{N_2}\frac{N_2}{ N_1}\right)\ge 
N_2 \ell \left(\frac{|X' (t)|}{N_2}\right).
\end{equation}
(A1,iii) implies (\ref{3.12.0427}) (see  (\ref{2.5.0610})).

We prove
\begin{equation}\label{3.14.0427}
T(P_0,P_1)\le \tilde V_2(P_0,P_1),
\end{equation}
which completes the proof of (i) from Theorem \ref{thm2.1.0416}.
The following implies (\ref{3.14.0427}): for $X(\cdot)\in \mathcal{A}_\infty(P_0,P_1)$,
from (A2,i),
\begin{equation}\label{3.15.0427}
N_2\int_0^1\ell \left(\frac{|X' (t)|}{N_2}\right)dt \ge \ell\left(|X(1)-X(0)|\right),
\end{equation}
in the same way as (\ref{3.2.0416}),
where $N_2=N_2(X')$.
Indeed, if $||X'||_1 >0$, then
\begin{equation}\label{3.16.0427}
\ell \left(\frac{|X' (t)|}{N_2}\right)
=\ell\left(\frac{|X' (t)|}{N_2||X'||_1 }||X'||_1 \right)
\ge \frac{|X' (t)|}{N_2||X'||_1 }\ell\left(||X'||_1 \right).
\end{equation}

We prove (ii).
For $X(\cdot)\in \mathcal{A}_{1,opt}(P_0,P_1)$, the equality holds in (\ref{3.13.0428}) from  
(\ref{2.6.0427}),
which implies that 
$X(\cdot)\in \mathcal{A}_{2,opt}(P_0,P_1)$.

For $X(\cdot)\in \mathcal{A}_{2,opt}(P_0,P_1)$, the equalities hold in (\ref{3.15.0427})--(\ref{3.16.0427}),
since  from  (\ref{2.1.0415}) and (\ref{2.6.0427}),
$\tilde V_2(P_0,P_1)=T(P_0,P_1)$.
This implies that $N_2(X')=N_1(X')$ from (A2,ii).
In particular, $X(\cdot)\in \mathcal{A}_{1,opt}(P_0,P_1)$  from (\ref{2.6.0427}).
\end{proof}

We prove Proposition  \ref{pp2.1}.
\begin{proof}
(Proposition  \ref{pp2.1})
If $|X_1-X_0|\ge r_0$, then take a random variable $Y$ such that the following holds:
$$C:=|Y-X_1|=|Y-X_0|=1+\frac{|X_1-X_0|}{2}.$$
\begin{equation}
Y(t):=
\begin{cases}
\displaystyle X_0+2t(Y-X_0),&0\le t\le \frac{1}{2},\\
\displaystyle Y+(2t-1)(X_1-Y),&\frac{1}{2}\le t\le 1.
\end{cases}
\end{equation}
Then $Y(t)=X_t, t=0,1$ and the following holds under our assumption:
\begin{equation}\label{3.18.0913}
N_2(Y')\int_0^1 \ell\left(\frac{|Y'(t)|}{N_2(Y')}\right)dt=\ell(2C)<\ell(|X_1-X_0| ),
\end{equation}
since  $|Y'(t)|=||Y'||_\infty =||Y'||_1=2C$.

If $|X_1-X_0|< r_0$, then 
\begin{equation}
Y(t):=X_0+t(X_1-X_0),\quad 0\le t\le  1.
\end{equation}
Then $Y(t)=X_t, t=0,1$ and the following holds:
\begin{equation}\label{3.20.0913}
N_2(Y')\int_0^1 \ell\left(\frac{|Y'(t)|}{N_2(Y')}\right)dt=\ell(|X_1-X_0| ).
\end{equation}
since $|Y'(t)|=||Y'||_\infty =||Y'||_1=|X_1-X_0|$.

From (\ref{3.18.0913}) and (\ref{3.20.0913}), under our assumption,
the following holds:
\begin{eqnarray}
\tilde V_2(P_0,P_1)
&\le &E\left [N_2(Y')\int_0^1 \ell\left(\frac{|Y'(t)|}{N_2(Y')}\right)dt\right]\\
&<&E[\ell(|X_1-X_0| )]=T(P_0,P_1).\notag
\end{eqnarray}
\end{proof}

We  prove Theorem \ref{thm2.1}.
\begin{proof}
(Theorem \ref{thm2.1})
We first prove (i).
We prove
\begin{equation}\label{3.1}
T^V(P_0,P_1;B)\le V(P_0,P_1;B).
\end{equation}
If $(X(\cdot),M)\in \mathcal{A}_\infty(P_0,P_1;B)$, then
$P^M\in B$ and 
\begin{eqnarray}
|X(1)-X(0)|&\le &\int_0^1|X'(t)|dt\le M,\label{3.2}\\
\int_0^1\ell\left(\frac{|X'(t)|}{M}M\right)dt
&\ge &\int_0^1\frac{|X'(t)|}{M}\ell (M)dt\label{3.14.0425}\\
&\ge &\frac{\ell (M)}{M}|X (1)-X(0)|,\notag
\end{eqnarray}
from (A1,i), provided $M>0$, which implies (\ref{3.1}).

We prove  
\begin{equation}\label{3.3}
V(P_0,P_1;B)\le T^V(P_0,P_1;B).
\end{equation}
If $(X_0, X_1,M)\in \Pi_\infty (P_0, P_1;B)$, then
$(X_M(\cdot):=X(\cdot;X_0, X_1,A_M),M)\in \mathcal{A}_\infty(P_0,P_1;B)$,
where $A_M:=[0,|X_1-X_0|/M]$ if $M>0$ and $=\{0\}$ if $M=0$.
In the case where $M\ne 0$,
$$|X_M'(t)|=\begin{cases} 
M, &\displaystyle 0<t<\frac{|X_1-X_0|}{M},\\
0,&\displaystyle \frac{|X_1-X_0|}{M}< t< 1,
\end{cases}
$$
\begin{equation}\label{2.8-1}
\int_0^1\ell (|X_M'(t)|)dt
=\int_0^{\frac{|X_1-X_0|}{M}}\ell (M)dt= \frac{\ell (M)}{M}|X_1-X_0|,
\end{equation}
which implies (\ref{3.3}).
(\ref{3.1}) and (\ref{3.3}) imply (\ref{2.1}).

We prove (ii).
Since $(X_0,X_1,M)\in \Pi_{T^V,opt}(P_0,P_1;B)$,
the following holds: from 
(\ref{2.1}),
\begin{eqnarray}
V(P_0,P_1;B)&\le& E\left[\int_0^1\ell(|X'(t;X_0,X_1,A)|)dt\right]\\
&=&E\left[\frac{\ell (M)}{M}|X_1-X_0|;M>0\right]\notag\\
&=&T^V(P_0,P_1;B)=V(P_0,P_1;B).\notag
\end{eqnarray}
Indeed, if $M\ge |X_1- X_0|>0$, then $|A|>0$ and 
$$|X'(t;X_0,X_1,A)|=\frac{|X_1-X_0|}{|A|}=M\quad \hbox{ on } A, \quad dt{\rm -a.e.}.$$

We prove  the first part of (iii).
Since $(X(\cdot),M)\in \mathcal{A}_{opt}(P_0,P_1;B)$,
the following holds: from  (\ref{2.1}) and (\ref{3.2})--(\ref{3.14.0425}),
\begin{eqnarray}\label{3.5}
T^V(P_0,P_1;B)
&=&V(P_0,P_1;B)\\
&=& E\left[\int_0^1\ell(|X'(t)|)dt\right]\notag\\
&\ge& E\left[\frac{\ell (M)}{M}|X (1)-X(0)|; M>0\right]
\notag\\
&\ge& T^V(P_0,P_1;B).\notag
\end{eqnarray}

We prove  the second part of (iii).
For  $(X(\cdot),M)\in \mathcal{A}_{opt}(P_0,P_1;B)$, from (A1,ii),
\begin{equation}\label{3.6}
|X' (t)|=0  \hbox{ or }M, \quad dtdP{\rm -a.e.,}
\end{equation}
\begin{equation}\label{2.10}
|X (1)-X(0)|=||X'||_1
,\quad {\rm a.s.},
\end{equation}
since the equality holds in (\ref{3.14.0425}) from (\ref{2.1}).
The following can be proved in the same way as (\ref{3.9.0416}):
\begin{equation}\label{2.19}
X' (t)=\frac{X (1)-X(0)}{|(X')^{-1}(\mathbb{R}^d\backslash\{0\})|}, \quad \hbox{ on }(X')^{-1}(\mathbb{R}^d\backslash\{0\}),  \quad dtdP{\rm - a.e.}.
\end{equation}
Indeed, replace $N_1|X (1)-X(0)|$ by $M$ in  (\ref{3.8.0416}).
(\ref{3.11.0416}) also implies the following:
\begin{equation}
|X (1)-X(0)|=M\times |(X')^{-1}(\mathbb{R}^d\backslash\{0\})|,\quad {\rm a.s.},
\end{equation}
which completes the proof.
\end{proof}

We  prove Corollary \ref{thm2.5.0628}.
\begin{proof}
(Corollary \ref{thm2.5.0628})
From Theorem \ref{thm2.1}, (\ref{2.17.0628}) can be obtained by the following:
\begin{equation}\label{3.27.0715}
T^V(P_0,P_1;\mathcal{P}([0,\infty)))=C_\ell \cdot T_1(P_0,P_1),
\end{equation} 
\begin{equation}\label{3.28.0715}
V(P_0,P_1;\mathcal{P}([0,\infty)))=\inf\left\{E\left[\int_0^1 \ell (|X'(t)|)dt\right]:
X(\cdot) \in  \mathcal{A}_\infty(P_0, P_1)\right\}.
\end{equation} 
(\ref{3.27.0715})  can be proved by the following.
If $(X_0,X_1,M)\in \Pi_{\infty}(P_0,P_1;\mathcal{P}([0,\infty)))$,
then $(X_0,X_1)\in \Pi(P_0,P_1)$ and 
$$E[\ell (M)M^{-1}|X_1-X_0|;M>0]\ge C_\ell E[|X_1-X_0|;M>0]= C_\ell E[ |X_1-X_0|].$$
If $(X_0,X_1)\in \Pi(P_0,P_1)$, then for $R>0$,
$(X_0,X_1,\max (|X_1-X_0|,R))\in \Pi_{\infty}(P_0,P_1;\mathcal{P}([0,\infty)))$,
and by the dominated convergence theorem,
$$E\left[\frac{\ell (\max (|X_1-X_0|,R))}{\max (|X_1-X_0|,R)}|X_1-X_0|\right]\to  C_\ell E[ |X_1-X_0|],
\quad R\to\infty.$$
(\ref{3.28.0715})  can be proved by the following.
If $(X(\cdot),M)\in \mathcal{A}_\infty(P_0,P_1;\mathcal{P}([0,\infty)))$,
then $X(\cdot) \in  \mathcal{A}_\infty(P_0, P_1)$.
If $X(\cdot) \in  \mathcal{A}_\infty(P_0, P_1)$, then  $(X(\cdot),||X'||_\infty)\in \mathcal{A}_\infty(P_0,P_1;\mathcal{P}([0,\infty)))$.

For $X(\cdot) \in  \mathcal{A}_\infty(P_0, P_1)$, from (A1,i),
\begin{equation}
\int_0^1 \ell (|X'(t)|)dt\ge C_\ell \int_0^1 |X'(t)|dt,
\end{equation} 
where the equality holds if and only if 
$$\ell (|X'(t)|)=C_\ell |X'(t)|,\quad dtdP-a.e..$$
If $P_0\ne P_1$, then $P(||X'||_\infty=0)<1$, which implies (\ref{2.13.0718}).
\end{proof}

\section{Appendix}

In this section, we state the proofs of (\ref{1.6}),  (\ref{1.9.0416}), and (\ref{1.11.0508}).

We prove (\ref{1.6}).
For $(X_0,X_1)\in \Pi(P_0, P_1)$ such that $E[L(X_1-X_0)]$ is finite,
\begin{equation*}
Y_n(t):=
\begin{cases}
\displaystyle X_0+n(X_1-X_0)t,&\displaystyle 0\le t\le \frac{1}{n},\\
\displaystyle X_1,&\displaystyle \frac{1}{n}\le t\le 1.
\end{cases}
\end{equation*}
Then  $Y_n\in \mathcal{A}(P_0, P_1)$, and
\begin{eqnarray*}
0\le V(P_0,P_1)\le E\left[\int_0^1 L(Y_n'(t))dt\right]&=&E[n^{-1}L(n(X_1-X_0))]\\
&\to& 0, \quad n\to\infty\notag
\end{eqnarray*}
by Lebesgue's dominated convergence theorem since 
$$n^{-1}L(n(X_1-X_0))\le L(X_1-X_0).$$

We prove (\ref{1.9.0416}).
For $(X_0,X_1) \in\Pi_{T,opt}(P_0,P_1)$,
$X(\cdot)$ defined by (\ref{3.6.0416}) is a minimizer of $V(P_0,P_1)$ by Jensen's inequality and 
$N_1(X')=1$, provided $V(P_0,P_1)$ is finite.
The following implies that  (\ref{1.9.0416}) holds: for $N\ge 1$,
$$L(x)=L\left(\frac{1}{N}Nx+(1-\frac{1}{N})0\right)\le \frac{1}{N}L(Nx)+(1-\frac{1}{N})L(0)=\frac{1}{N}L(Nx)$$
since $L:\mathbb{R}^d\rightarrow [0,\infty)$ is convex and $L(0)=0$.

We prove (\ref{1.11.0508}).
For $T\ge 1, \varphi\in AC([0,T])$ such that $N_i(\varphi')_T=T$,
$$\varphi (T \cdot)\in AC([0,1]),\quad 
N_i(\varphi' (T \cdot))=N_i(\varphi')_T=T,
$$
\begin{eqnarray*}
\int_0^T\ell \left(|\varphi'(t)|\right)dt
&=&\int_0^1\ell \left(|\varphi'(Ts)|\right)Tds\\
&=&\int_0^1N_i(\varphi' (T \cdot))\ell \left(\frac{1}{N_i(\varphi' (T \cdot))}\left|\frac{d}{ds}\varphi(Ts)\right|\right)ds,
\end{eqnarray*}
which implies that  (l. h. s.) $\le$ (r. h. s.) in (\ref{1.11.0508}).

For $\varphi\in AC([0,1]), T\ge 1$,
$$\varphi \left(\frac{\cdot}{T}\right)\in AC([0,T]),\quad 
N_i\left(\varphi'\left(\frac{\cdot}{T}\right)\right)_T=N_i(\varphi'),
$$
\begin{eqnarray*}
\int_0^1N_i(\varphi')\ell \left(\frac{|\varphi'(t)|}{N_i(\varphi')}\right)dt
&=&\int_0^TN_i(\varphi')\ell \left(\frac{1}{N_i(\varphi')}\left|\varphi '\left(\frac{s}{T}\right)\right|\right)\frac{1}{T}ds\\
&=&\int_0^T\frac{N_i\left(\varphi'\left(\frac{\cdot}{T}\right)\right)_T}{T}\ell \left(\frac{T}{N_i\left(\varphi'\left(\frac{\cdot}{T}\right)\right)_T}\left|\frac{d}{ds}\varphi \left(\frac{s}{T}\right)\right|\right)ds\\
&=&\int_0^T\ell \left(\left|\frac{d}{ds}\varphi \left(\frac{s}{T}\right)\right|\right)ds,
\end{eqnarray*}
provided $T=N_i(\varphi')$.
This  implies that  (l. h. s.) $\ge$ (r. h. s.) in (\ref{1.11.0508}).

\end{document}